\documentclass[12pt]{amsart}
\usepackage{graphicx}
\usepackage{longtable}
\newtheorem{thm}{Theorem}

\newtheorem{defn}[thm]{Definition}
\newtheorem*{remark}{Remark}
\newtheorem*{ack}{Acknowledgements}

\newcommand{\Z}{\mathbb{Z}}
\newcommand{\Q}{\mathbb{Q}}
\newcommand{\F}{\mathbb{F}}

\renewcommand{\P}{\mathbb{P}}
\newcommand{\sfp}{{\rm sfp}}
\newcommand{\ord}{{\rm ord}}
\newcommand{\rank}{{\rm rank}}

\DeclareFontFamily{U}{wncy}{}
\DeclareFontShape{U}{wncy}{m}{n}{<->wncyr10}{}
\DeclareSymbolFont{mcy}{U}{wncy}{m}{n}
\DeclareMathSymbol{\Sha}{\mathord}{mcy}{"58}
\begin{document}

\title{Three consecutive almost squares}

\author{Jeremy Rouse}
\address{Department of Mathematics, Wake Forest University, Winston-Salem, North Carolina 27109 USA}
\email{rouseja@wfu.edu}
\author{Yilin Yang}
\address{Department of Mathematics, University of Michigan, Ann Arbor,
Michigan 48109 USA}
\email{yangy210v@gmail.com}

\begin{abstract}
Given a positive integer $n$, we let $\sfp(n)$ denote the squarefree part
of $n$. We determine all positive integers $n$ for which
$\max \{ \sfp(n), \sfp(n+1), \sfp(n+2) \} \leq 150$ by relating the problem
to finding integral points on elliptic curves. We also prove that
there are infinitely many $n$ for which
\[
  \max \{ \sfp(n), \sfp(n+1), \sfp(n+2) \} < n^{1/3}.
\]
\end{abstract}

\maketitle

\section{Introduction}

The positive integers $48$, $49$ and $50$ are consecutive, and are
``almost squares'', namely $48 = 3 \cdot 4^{2}$, $49 = 7^{2}$ and
$50 = 2 \cdot 5^{2}$. Does this phenomenon ever occur again? That is,
is there a positive integer $n > 48$ for which
\begin{align*}
  n &= 3x^{2}\\
  n+1 &= y^{2}\\
  n+2 &= 2z^{2}
\end{align*}
has an integer solution? The answer is no. One perspective on this is given
by Cohen in Section 12.8.2 of \cite{Cohen}, where this problem is solved
using linear forms in logarithms. Another approach is to recognize that
the pair of equations $2z^{2} - y^{2} = y^{2} - 3x^{2} = 1$ define
an intersection of two quadrics in $\P^{3}$ and hence define a curve of genus
$1$. Siegel proved that there are only finitely many integer points on
a genus $1$ curve, and hence the system of equations above has only
finitely many solutions. This method of solving simultaneous Pell
equations is considered in \cite{Tzanakis}.

We consider more generally the problem of, given integers $a$, $b$ and
$c$, finding integers $n$ for which $n = ax^{2}$, $n+1 = by^{2}$ and
$n+2 = cz^{2}$.  Multiplying these equations gives $n(n+1)(n+2) =
abc(xyz)^{2}$ and is related to the problem of finding consecutive
integers that multiply to an ``almost'' perfect square. There is an
extensive literature on related problems. On one hand, Erd\H{o}s and
Selfridge proved in \cite{ErdSelf} that no product of consecutive
integers is a perfect power (a result generalized to arithmetic
progressions by others, see \cite{Ben1} and \cite{GHP}). In another
direction, Gy\H{o}ry proved \cite{Gyory} that the equation
$n(n+1)\cdots(n+k-1) = bx^{l}$ has no integer solutions with $n > 0$
if the greatest prime factor of $b$ does not exceed $k$.

For a positive integer $n$, define $\sfp(n)$ to be the ``squarefree
part'' of $n$, the smallest positive integer $a$ so that $a | n$ and
$\frac{n}{a}$ is a perfect square. In this paper, we consider positive
integers $n$ for which $\sfp(n)$, $\sfp(n+1)$ and $\sfp(n+2)$ are all
small. Our first result is a classification of values of $n$
for which $\sfp(n)$, $\sfp(n+1)$ and $\sfp(n+2)$ are all $\leq 150$.
We say that such an $n$ is ``non-trivial'' if
$\sfp(n) < n$, $\sfp(n+1) < n+1$ and $\sfp(n+2) < n+2$.

\begin{thm}
\label{upto150}
There are exactly $25$ non-trivial $n$ for which $\sfp(n) \leq 150$, $\sfp(n+1) \leq 150$ and $\sfp(n+2) \leq 150$, the largest of which is $n = 9841094$. 
\end{thm}
\begin{remark}
A table of all $25$ values of $n$ is given in Section~\ref{table}.
\end{remark}

After this manuscript was complete, the authors were informed that the
essence of this result follows from work of Michael Bennett and Gary
Walsh (see \cite{BenWal}). Given positive integers $a$, $b$ and $c$, a
positive integer solution to $n = ax^{2}$, $n+1 = by^{2}$, $n+2 =
cz^{2}$ yields $(n+1)^{2} - n(n+2) = b^{2} y^{4} - (ac) (xz)^{2} =
1$. Bennett and Walsh prove that the system $b^{2} Y^{4} - d X^{2} =
1$ has at most one solution in positive integers. Moreover, if $T + U
\sqrt{d}$ is a fundamental solution of Pell's equation $X^{2} - dY^{2}
= 1$ define for $k \geq 1$, $T_{k} + U_{k} \sqrt{d} = (T + U
\sqrt{d})^{k}$. Then $T_{k} = bx^{2}$ for some integer $x \in \Z$ for
at most one $k$, and if any such $k$ exists, then $k$ is the smallest
positive integer for which $T_{k}$ is a multiple of $b$. Bennett and
Walsh's approach uses linear forms in logarithms, while our approach
relies heavily on the theory of elliptic curves. In particular, if $C$
is the curve defined by the two equations $by^{2} - ax^{2} = 1, cz^{2}
- by^{2} = 1$, then the Jacobian of $C$ is isomorphic to $E : y^{2} =
x^{3} - (abc)^{2} x$.

We rule out many of the $778688$
candidates for the triple $(a,b,c)$ by checking to see whether there
are integer solutions to each of the three equations $by^{2} - ax^{2} =
1$, $cz^{2} - b^{2} = 1$ and $cz^{2} - ax^{2} = 2$. We also test $C$
for local solvability, and use Tunnell's theorem (see
\cite{Tunnell}) to determine if the rank of $E$ is positive.  Finally,
we use the surprising property that the natural map from $C$ to $E$
sends an integral solution on $C$ to an integral point on $E$.  It
suffices therefore to compute all the integral points on $E$ (using
Magma \cite{Magma}), which
requires computing generators of the Mordell-Weil group. In many cases
this is straightforward, but a number of cases require more involved
methods (12-descent, computing the analytic rank, and the use of Heegner points). 

Given the existence of large solutions, it is natural to ask how large\\
$\max \{ \sfp(n), \sfp(n+1), \sfp(n+2) \}$ can be as a function of $n$.
This is the subject of our next result.
\begin{thm}
\label{family}
There are infinitely many positive integers $n$ for which
\[
  \max \{ \sfp(n), \sfp(n+1), \sfp(n+2) \} < n^{1/3}.
\]
\end{thm}
\begin{remark}
  The following heuristic suggests that the exponent $1/3$ above is
  optimal.  A partial summation argument shows that the number of
  positive integers $n \leq x$ with $\sfp(n) \leq z$ is $\frac{12
    \sqrt{xz}}{\pi^{2}} + O(\sqrt{x} \log(z))$. Assuming that $n$,
  $n+1$ and $n+2$ are ``random'' integers it follows that $\sfp(n),
  \sfp(n+1)$ and $\sfp(n+2)$ are all $\leq n^{\alpha}$ with
  probability about $n^{3 (\alpha - 1)/2}$. Therefore, the
  ``expected'' number of $n$ for which $\max \{ \sfp(n), \sfp(n+1),
  \sfp(n+2) \} \leq n^{\alpha}$ is infinite if $3 (\alpha - 1)/2 \geq
  -1$ and finite otherwise.
\end{remark}

\begin{ack}
This work arose out of the second author's undergraduate thesis
at Wake Forest University. The authors used Magma \cite{Magma} V2.20-9
for computations. The authors would also like to thank the anonymous referee
for a number of very helpful comments.
\end{ack}

\section{Background}

We denote by $\mathbb{Q}_{p}$ the field of $p$-adic numbers.  A
necessary condition for a curve $X/\mathbb{Q}$ to have a rational
solution is for it to have such a solution in $\mathbb{Q}_{p}$ for all
primes $p$.

If $d$ is an integer which is not a perfect square, let
$\Z[\sqrt{d}] = \{ a + b \sqrt{d} : a, b \in \Z \}$.  This is a (not
necessarily maximal) order in the quadratic field $\Q[\sqrt{d}]$. Let
$N : \Z[\sqrt{d}] \to \Z$ be the norm map given by
$N(a + b \sqrt{d}) = a^{2} - db^{2}$.

We now describe some background about elliptic curves. For our
purposes an elliptic curve is a curve of the shape
\[
  E : y^{2} = x^{3} + ax + b
\]
where $a, b \in \Q$. Let $E(\Q)$ be the set of pairs $(x,y)$ of rationals
numbers that solve the equation, together with the ``point at infinity''. This
set has a binary operation on it: given $P, Q$ in $E(\Q)$, the line $L$
through $P$ and $Q$ intersects $E$ in a third point $R = (x,y)$. The point
$P+Q$ is defined to be $(x,-y)$. This binary operation endows $E(\Q)$ with
the structure of an abelian group.

The Mordell-Weil theorem (see \cite{Silverman}, Theorem VIII.4.1, for example)
states the following.
\begin{thm}\label{TZ}
For any elliptic curve $E/\Q$, the group $E(\Q)$ is finitely generated. More,
precisely,
\[
E( \mathbb{Q} ) \cong E(\mathbb{Q})_{\rm tors} \times \mathbb{Z}^{\rank(E(\Q))},
\]
where $E(\mathbb{Q})_{tors}$ is the (finite) torsion subgroup.
\end{thm}

There is in general no algorithm which is proven to compute the
rank of $E$, (see Section~3 of
Rubin and Silverberg's paper \cite{Rubin}) but there are a number of
procedures which work well in practice for relatively simple curves
$E$. We will reduce the problem of solving $n = ax^{2}$, $n+1 =
by^{2}$ and $n+2 = cz^{2}$ to finding points $(X,Y)$ on the curve
\[
  E : Y^{2} = X^{3} - (abc)^{2} X,
\]
with $X, Y \in \Z$. A theorem of Siegel (see Theorem IX.4.3 of
\cite{Silverman}) states that there are only finitely many points in
$E(\Q)$ with both coordinates integral. There are effective and
practical algorithms (see \cite{ST}, \cite{IntegralPoints}, \cite{Sm} and
\cite{Tz}) to determine the set
of integral points, provided the rank $r$ can be computed and a set of
generators for $E(\Q)$ found. Given a point $P = (x,y) \in E(\Q)$, the
``naive height'' of $P$ is defined by writing $x = \frac{a}{b}$ with
$a, b \in \Z$ with $\gcd(a,b) = 1$ and defining $h(P) = \log \max \{
|a|, |b| \}$. The ``canonical height'' of $P$ is defined to be
\[
  \hat{h}(P) = \lim_{n \to \infty} \frac{h(2^{n} P)}{4^{n}}.
\]

The $L$-function of $E/\Q$ is defined to be
\[
  L(E,s) = \sum_{n=1}^{\infty} \frac{a_{n}(E)}{n^{s}}
  = \prod_{p} \left(1 - a_{p}(E) p^{-s} + \varepsilon(p) p^{1-2s}\right)^{-1}
\]
where $a_{p}(E) = p + 1 - \# E(\F_{p})$, and 
\[
  \varepsilon(p) = \begin{cases}
    1 & \text{ if } p \text{ does not divide the conductor of } E,\\
    0 & \text{ otherwise. }
\end{cases}
\]

The Birch and Swinnerton-Dyer conjecture states
that $\ord_{s=1} L(E,s) = \rank(E(\Q))$, and moreover that
\[
  \lim_{s \to 1} \frac{L(E,s)}{(s-1)^{r}} = \frac{\Omega R(E/\Q) \Sha(E/\Q) \prod_{p} c_{p}}{|T|^{2}}.
\]
In our case, $\Omega$ is twice the real period, $R(E/\Q)$ is the regulator
of $E(\Q)$ computed using the function $\hat{h}$ above, the $c_{p}$ are the
Tamagawa numbers, and $\Sha(E/\Q)$ is the Shafarevich-Tate group.
The completed $L$-function $\Lambda(E,s) = N^{s/2} (2 \pi)^{-s} \Gamma(s) L(E,s)$
satisfies the function equation $\Lambda(E,s) = w_{E} \Lambda(E,2-s)$,
where $w_{E}$ is the root number of $E$. Note that
$w_{E} = 1$ implies that $\ord_{s=1} L(E,s)$ is even, and $w_{E} = -1$
if $\ord_{s=1} L(E,s)$ is odd.

The best partial result in the direction of the Birch and Swinnerton-Dyer
conjecture is the following.

\begin{thm}[Gross-Zagier \cite{GrossZagier}, Kolyvagin \cite{Kolyvagin}, et al.]
\label{bsdpartial}
Suppose that $E/\Q$ is an elliptic curve and $\ord_{s=1} L(E,s) = 0$ or $1$.
Then, $\ord_{s=1} L(E,s) = \rank(E(\Q))$.
\end{thm}
The work of Bump-Friedberg-Hoffstein \cite{BFH} or Murty-Murty \cite{MM}
is necessary to remove a condition imposed in the work of Gross-Zagier
and Kolyvagin. 

\section{Proof of Theorem~\ref{family}}

\begin{proof}
Motivated by the observation that $8388223 = 127 \cdot 257^{2}$
and $8388225 = 129 \cdot 255^{2}$, we find a parametric family
of solutions where $n = a \cdot b^{2}$ and $n+2 = (a+2) \cdot (b-2)^2$.

If we write $(4 + \sqrt{13}) (649+180 \sqrt{13})^{m} = x_{m} + y_{m} \sqrt{13}$
where $x_{m}$ and $y_{m}$ are integers, then $x_{m}^{2} - 13 y_{m}^{2} = 3$ for
all $m \geq 0$. We have that $x_{0} = 4$,
$x_{1} = 4936$, $y_{0} = 1$, $y_{1} = 1369$ and
\[
  x_{m} = 1298 x_{m-1} - x_{m-2}, y_{m} = 1298 y_{m-1} - y_{m-2}.
\]
It is easy to see that $x_{m}$ is periodic modulo $32$ with period $8$
and from this it follows that $x_{8m+7} \equiv 0 \pmod{32}$ for all $m \geq 0$. 
Set $a_{m} = x_{8m+7}/2$ and $n = 4 a_{m}^{3} - 3 a_{m} - 1$. Then we have
\begin{align*}
  n &= (a_{m} - 1) (2 a_{m} + 1)^{2},\\
  n+1 &= a_{m} (4 a_{m}^{2} - 3) = 13 a_{m} y_{8m+7}^{2}, \text{ and }\\
  n+2 &= (a_{m} + 1) (2 a_{m} - 1)^{2}.\\
\end{align*}
Since $16 | a_{m}$, $\max \{ \sfp(n), \sfp(n+1), \sfp(n+2) \} \leq \max \{ a_{m} - 1, \frac{13 a_{m}}{16}, a_{m} + 1 \} = a_{m} + 1 < n^{1/3}$.
\end{proof}
\begin{remark}
The polynomial $4x^{3} - 3x$ used in the proof above is the
Chebyshev polynomial $T_{3}(x)$. This explains why $T_{3}(x) - 1$ and $T_{3}(x) + 1$ both have a double zero.
\end{remark}

\section{Proof of Theorem~\ref{upto150}}

Since there are $92$ squarefree integers $\leq 150$, 
there are $778688 = 92^{3}$ possibilities for the triple $(a,b,c)$.
We first test four things before searching for integral
points on $E : y^{2} = x^{3} - (abc)^{2} x$. Suppose that $n = ax^{2}$,
$n+1 = by^{2}$ and $n+2 = cz^{2}$ is an integral solution to the system of 
equations 
\begin{align}
\label{eq1}  by^{2} - ax^{2} &= 1,\\ 
\label{eq2}  cz^{2} - by^{2} &= 1,\\
\label{eq3}  cz^{2} - ax^{2} &= 2.
\end{align}

\subsection{Greatest common divisor conditions}

If $(x,y,z)$ is an integer solution to \eqref{eq1}, \eqref{eq2}
and \eqref{eq3}, then $\gcd(a,b) = 1$, $\gcd(b,c) = 1$ and
$\gcd(a,c) = 1 \text{ or } 2$. This reduces the number of triples to
consider to $425639$.

\subsection{Norm equations}

We next check that each of the equations $by^{2} - ax^{2} = 1$,
$cz^{2} - by^{2} = 1$, and $cz^{2} - ax^{2} = 2$ has an integer solution.
The equation $by^{2} - ax^{2} = 1$ has an integer solution if and only if
$Y^{2} - abx^{2} = b$ has an integer solution. This equation has a solution
if and only if $\Z[\sqrt{ab}]$ has an element of norm $b$. We use Magma's
routine {\tt NormEquation} to test this. After these tests have been
made, there are $2188$ possibilities for $(a,b,c)$ that remain.

\subsection{Local solvability of $C$}

Let $C \subseteq \P^{3}$ be the curve defined by the two equations
$by^{2} - ax^{2} = w^{2}$, $cz^{2} - by^{2} = w^{2}$. This curve must
have a rational point on it in order for there to be a non-trivial
solution. We check whether $C$ has points $(x : y : z : w)$ in $\Q_{p}$
for all primes $p$ dividing $2abc$. This is done via Magma's {\tt
  IsLocallySolvable} routine. This eliminates $244$ possibilities.

\subsection{Rank of the elliptic curve $E$}

Let $E : y^{2} = x^{3} - (abc)^{2} x$. Define a map $M$ from
non-trivial solutions to $C$, represented in the form
$n = ax^{2}$, $n+1 = by^{2}$ and $n+2 = cz^{2}$ to $E$, given by
\[
  M(x,y,z,n) = ((n+1)(abc), (abc)^{2} xyz).
\]
If $n > 0$, then $xyz > 0$ and so $M(x,y,z,n) \in E(\Q)$ is an integral
point, and one with a non-zero $y$-coordinate. We now use the following
results about the family of elliptic curves $E$ we consider.

\begin{defn}
A natural number $N$ is called \emph{congruent} if there exists a right 
triangle with all three sides rational and area $N$.
\end{defn}
Consider the elliptic curve $E$ over $\mathbb{Q}$ given by:
\[
E(\mathbb{Q}): y^{2}=x^{3}-N^{2} x.
\]
We have the following result.
\begin{thm}[Proposition I.9.18 of \cite{Koblitz}]
\label{congthm}
The number $N$ is congruent if and only if the rank of $E$ is positive.
\end{thm}
Proposition I.9.17 of \cite{Koblitz} implies that the only points of
finite order on $E$ are $(0,0)$, $(\pm N, 0)$ and the point at infinity.
If for some $n > 0$ we have $n = ax^{2}$,
$n+1 = by^{2}$ and $n+2 = cz^{2}$, then the point $M(x,y,z,n)$ on
$E : Y^{2} = X^{3} - (abc)^{2} X$ is, according to the above, a non-torsion
point, and hence $E$ has positive rank, and consequently, by 
Theorem~\ref{congthm}, $abc$ is a congruent number.

\begin{thm} [Tunnell 1983] If $N$ is squarefree and odd, and
\[
\# \left\{ x,y,z \in \mathbb{Z} \vert N =2x^{2}+y^{2}+32z^2 \right\} \ne
\frac{1}{2} \# \left\{ x,y,z \in \mathbb{Z} \vert
  N=2x^{2}+y^{2}+8z^2 \right\},
\]
then $N$ is not congruent. If $N$ is squarefree and even, and
\[
\# \left\{ x,y,z \in \mathbb{Z} : \frac{N}{2}=4x^{2}+y^{2}+32z^2 \right\} \ne
\frac{1}{2}  \# \left\{ x,y,z \in \mathbb{Z} : \frac{N}{2}=4x^{2}+y^{2}+8z^2 \right\} ,
\]
then
$N$ is not congruent.
\end{thm}
To use Tunnell's theorem, we compute the generating function for the
number of representations of $N$ by $2x^{2} + y^{2} + 8z^{2}$ as
\[
  \sum_{x, y, z \in \Z} q^{2x^{2} + y^{2} + 8z^{2}}
  = \left(\sum_{x \in \Z} q^{2x^{2}}\right) \left(\sum_{y \in \Z} q^{y^{2}}\right) \left(\sum_{z \in \Z} q^{8z^{2}}\right),
\]
as well as the representations of $N$ by $2x^{2}+y^{2}+32z^{2}$, $4x^{2}+y^{2}+8z^{2}$ and $4x^{2}+y^{2}+32z^{2}$. Actually, for our purposes, it suffices
to compute
\[
  \sum_{x=-X}^{X} \sum_{y=-Y}^{Y} \sum_{z=-Z}^{Z}
  q^{2x^{2} + y^{2} + 8z^{2}},
\]
where $X = \lfloor \sqrt{150^{3}/2} \rfloor$, $Y = \lfloor
\sqrt{150^{3}} \rfloor$, and $Z = \lfloor \sqrt{150^{3}/8} \rfloor$.

The use of Tunnell's theorem rules out $530$ of the $1944$ remaining cases.
In the case that the hypothesis of Tunnell's theorem is false, the Birch
and Swinnerton-Dyer conjecture predicts that $E(\Q)$ does have positive
rank. In this case, we proceed to the next step.

\subsection{Computing integral points}

Once the program determines that the elliptic curve $E(\mathbb{Q}):
Y^{2} = X^{3} - (abc)^{2} X$ has positive rank, we let Magma
attempt to compute the integral points on the curve, using the
routine {\tt IntegralPoints}, which is based on the method developed
in \cite{IntegralPoints} and \cite{ST}. If this routine does not succeed
within $15$ minutes, which is sufficient time to perform a 4-descent,
we abort the computation.

In $1377$ of the $1414$ cases that remain, Magma is able to determine
the Mordell-Weil group and determine all of the integral points within
$15$ minutes. Once the integral points are determined, we check to see
if they are in the image of the map $M(x,y,z,n)$ and if so, whether
they correspond to a non-trivial solution.

There are $37$ more difficult cases that remain, and in each
case we are able to use other methods to determine the Mordell-Weil group.
A table of these cases and the generators of the Mordell-Weil group
are given at the page {\tt http://users.wfu.edu/rouseja/MWgens.html}. 

Of the $37$ cases, there are four cases with root number $-1$ and
rank $\leq 3$ for which one point of infinite order has low height.
In these, we numerically compute $L'(E,1)$ and show that it is
nonzero. Theorem~\ref{bsdpartial} proves the rank is one in these
cases. One of these cases is $a =
139$, $b = 89$ and $c = 109$. In this case, $E$ has conductor $\approx
5.82 \cdot 10^{13}$ and computing $L'(E,1)$ takes about $5$ hours.

In the remaining $33$ cases, rather than use the $8$-descent (the
default if the {\tt IntegralPoints} command were to keep running), we
use the 12-descent routines in Magma (due to Tom Fisher \cite{Fisher})
to find points of large height. The curve with
$a = 137$, $b = 109$ and $c = 101$ has a generator with canonical height
$1234$, which is found after searching on the $12$-covers for about $90$
minutes. This approach is successful in all but one case.

\subsection{Heegner points}

All but one of the $778688$ original cases are handled by the methods
of the previous sections. The remaining case is $a = 67$, $b = 131$
and $c = 109$. The curve $E$ has root number $-1$, rank $\leq 1$ and
conductor $\approx 2.9 \cdot 10^{13}$. A long computation shows that
$L'(E,1) \approx 72.604$. This suggests, assuming the Birch and
Swinnerton-Dyer conjecture is true and $\Sha(E/\Q)$ is trivial, that a
generator of the Mordell-Weil group has canonical height about
$1692$. Searching for points on the 12-covers up to a height of
$3^{42} \cdot 10^{5} \approx 10^{25}$ does not succeed in finding
points.

For this final case, we use the method of Elkies described in
\cite{ElkiesHeegner}. This is a variant of the usual Heegner point
method and is quite fast on quadratic twists of curves with low
conductor. The modular curve $X_{0}(32)$ parametrizes pairs $(E,C)$,
where $E$ is an elliptic curve, and $C \subseteq E$ is a cyclic
subgroup of order $32$. It is well-known that $X_{0}(32)$ is
isomorphic to $y^{2} = x^{3} + 4x$, which in turn has a degree $2$ map
to $y^{2} = x^{3} - x$. Finding a rational point on
$y^{2} = x^{3} - D^{2} x$ is equivalent to finding a point on
$y^{2} = x^{3} - x$ with $x$ and $y/\sqrt{-D}$ both rational. If $E$
is an elliptic curve with complex multiplication whose endomorphism
ring $\mathcal{O}$ contains an ideal $I$ so that
$\mathcal{O}/I \cong \Z/32 \Z$, this naturally gives rise to a point
on $X_{0}(32)$, which is defined in the ring class field of
$\mathcal{O}$.  Taking the trace of such a point (in the Mordell-Weil
group) gives rise to a point on $X_{0}(32)$ over an imaginary
quadratic field.

In \cite{ElkiesHeegner}, Elkies gives a procedure for constructing rational
points on $y^{2} = x^{3} - D^{2} x$ with $D \equiv 7 \pmod{8}$. We must handle
one case with $D \equiv 5 \pmod{8}$. We let $\mathcal{O} = \Z[4 \sqrt{-D}]$.
This ring has an ideal $I = \langle 32, 4 + 4 \sqrt{-D} \rangle$ and
$\mathcal{O}/I \cong \Z/32 \Z$. We take a representative ideal $J$
for each element of the class group of $\mathcal{O}$. Thinking of $J$
as a lattice, $\mathbb{C}/IJ$ is an elliptic curve,
and $I/IJ$ is a subgroup of order $32$ on the curve. Thus,
$(\mathbb{C}/IJ, I/IJ)$ is a Heegner point on $X_{0}(32)$. Via the isomorphism
with $E : y^{2} = x^{3} + 4x$, we obtain a collection of points. Adding
these points together on $E$ gives a point defined in $\Q(\sqrt{-D})$,
and eventually a rational point on $y^{2} = x^{3} - D^{2} x$.
Running this computation in several simpler cases suggests that the resulting
point $P$ on $y^{2} = x^{3} - D^{2} x$ satisfies $mQ = P$ for some rational
point $Q$, where $m$ is equal to the number of divisors of $D$. Any
point $Q$ must have an $x$-coordinate of the form $x = em^{2}/n^{2}$,
where $e$ is a divisor of $2D$. We use this method to minimize the
amount of decimal precision needed to compute $P$.

For $D = 67 \cdot 131 \cdot 109$, the class number of $\mathcal{O}$ is
$3712$. Given that we expect a point with canonical height $1692$, we
compute the Heegner points and their trace in the manner described in
\cite{ElkiesHeegner} using $850$ digits of decimal precision. The
result is a point $Q$ with canonical height $1692.698$, and the
computation takes about 18 minutes. This point is a generator of the
Mordell-Weil group of $E : y^{2} = x^{3} - D^{2} x$. Once the
Mordell-Weil group is found, we check that there are no integral
points on $E$. This completes the proof of Theorem~\ref{upto150}.

\section{Table of $n$}
\label{table}

The following is a table of all $25$ positive integers $n$ with
$\sfp(k) < \max \{ k, 150 \}$ for $k = n$, $n+1$ and $n+2$.

\tiny
\begin{center}
\begin{tabular}{cccc}
$n$ & $\sfp(n)$ & $\sfp(n+1)$ & $\sfp(n+2)$\\
\hline
$48$ & $3$ & $1$ & $2$\\
$98$ & $2$ & $11$ & $1$\\
$124$ & $31$ & $5$ & $14$\\
$242$ & $2$ & $3$ & $61$\\
$243$ & $3$ & $61$ & $5$\\
$342$ & $38$ & $7$ & $86$\\
$350$ & $14$ & $39$ & $22$\\
$423$ & $47$ & $106$ & $17$\\
$475$ & $19$ & $119$ & $53$\\
$548$ & $137$ & $61$ & $22$\\
$845$ & $5$ & $94$ & $7$\\
$846$ & $94$ & $7$ & $53$\\
$1024$ & $1$ & $41$ & $114$\\
$1375$ & $55$ & $86$ & $17$\\
$1519$ & $31$ & $95$ & $1$\\
$1680$ & $105$ & $1$ & $2$\\
$3724$ & $19$ & $149$ & $46$\\
$9800$ & $2$ & $1$ & $58$\\
$31211$ & $59$ & $3$ & $13$\\
$32798$ & $62$ & $39$ & $82$\\
$118579$ & $19$ & $5$ & $141$\\
$629693$ & $53$ & $46$ & $55$\\
$1294298$ & $122$ & $19$ & $7$\\
$8388223$ & $127$ & $26$ & $129$\\
$9841094$ & $134$ & $55$ & $34$\\
\end{tabular}
\end{center}
\normalsize

\bibliographystyle{plain}
\bibliography{refs}

\def\cprime{$'$}
\begin{thebibliography}{10}

\bibitem{Ben1}
M.~A. Bennett, N.~Bruin, K.~Gy{\H{o}}ry, and L.~Hajdu.
\newblock Powers from products of consecutive terms in arithmetic progression.
\newblock {\em Proc. London Math. Soc. (3)}, 92(2):273--306, 2006.

\bibitem{BenWal}
Michael~A. Bennett and Gary Walsh.
\newblock The {D}iophantine equation {$b^2X^4-dY^2=1$}.
\newblock {\em Proc. Amer. Math. Soc.}, 127(12):3481--3491, 1999.

\bibitem{Magma}
Wieb Bosma, John Cannon, and Catherine Playoust.
\newblock The {M}agma algebra system. {I}. {T}he user language.
\newblock {\em J. Symbolic Comput.}, 24(3-4):235--265, 1997.
\newblock Computational algebra and number theory (London, 1993).

\bibitem{BFH}
Daniel Bump, Solomon Friedberg, and Jeffrey Hoffstein.
\newblock Nonvanishing theorems for {$L$}-functions of modular forms and their
  derivatives.
\newblock {\em Invent. Math.}, 102(3):543--618, 1990.

\bibitem{Cohen}
Henri Cohen.
\newblock {\em Number theory. {V}ol. {II}. {A}nalytic and modern tools}, volume
  240 of {\em Graduate Texts in Mathematics}.
\newblock Springer, New York, 2007.

\bibitem{ElkiesHeegner}
Noam~D. Elkies.
\newblock Heegner point computations.
\newblock In {\em Algorithmic number theory ({I}thaca, {NY}, 1994)}, volume 877
  of {\em Lecture Notes in Comput. Sci.}, pages 122--133. Springer, Berlin,
  1994.

\bibitem{ErdSelf}
P.~Erd{\H{o}}s and J.~L. Selfridge.
\newblock The product of consecutive integers is never a power.
\newblock {\em Illinois J. Math.}, 19:292--301, 1975.

\bibitem{Fisher}
Tom Fisher.
\newblock Finding rational points on elliptic curves using 6-descent and
  12-descent.
\newblock {\em J. Algebra}, 320(2):853--884, 2008.

\bibitem{IntegralPoints}
J.~Gebel, A.~Peth{\H{o}}, and H.~G. Zimmer.
\newblock Computing integral points on elliptic curves.
\newblock {\em Acta Arith.}, 68(2):171--192, 1994.

\bibitem{GrossZagier}
Benedict~H. Gross and Don~B. Zagier.
\newblock Heegner points and derivatives of {$L$}-series.
\newblock {\em Invent. Math.}, 84(2):225--320, 1986.

\bibitem{Gyory}
K.~Gy{\H{o}}ry.
\newblock On the {D}iophantine equation {$n(n+1)\cdots(n+k-1)=bx^l$}.
\newblock {\em Acta Arith.}, 83(1):87--92, 1998.

\bibitem{GHP}
K.~Gy{\H{o}}ry, L.~Hajdu, and {\'A}.~Pint{\'e}r.
\newblock Perfect powers from products of consecutive terms in arithmetic
  progression.
\newblock {\em Compos. Math.}, 145(4):845--864, 2009.

\bibitem{Koblitz}
Neal Koblitz.
\newblock {\em Introduction to elliptic curves and modular forms}, volume~97 of
  {\em Graduate Texts in Mathematics}.
\newblock Springer-Verlag, New York, second edition, 1993.

\bibitem{Kolyvagin}
Victor~A. Kolyvagin.
\newblock Finiteness of {$E({\bf Q})$} and {$\Sha(E,{\bf Q})$} for a subclass
  of {W}eil curves.
\newblock {\em Izv. Akad. Nauk SSSR Ser. Mat.}, 52(3):522--540, 670--671, 1988.

\bibitem{MM}
M.~Ram Murty and V.~Kumar Murty.
\newblock Mean values of derivatives of modular {$L$}-series.
\newblock {\em Ann. of Math. (2)}, 133(3):447--475, 1991.

\bibitem{Rubin}
Karl Rubin and Alice Silverberg.
\newblock Ranks of elliptic curves.
\newblock {\em Bull. Amer. Math. Soc. (N.S.)}, 39(4):455--474 (electronic),
  2002.

\bibitem{Silverman}
Joseph~H. Silverman.
\newblock {\em The arithmetic of elliptic curves}, volume 106 of {\em Graduate
  Texts in Mathematics}.
\newblock Springer, Dordrecht, second edition, 2009.

\bibitem{Sm}
Nigel~P. Smart.
\newblock {\em The algorithmic resolution of {D}iophantine equations},
  volume~41 of {\em London Mathematical Society Student Texts}.
\newblock Cambridge University Press, Cambridge, 1998.

\bibitem{ST}
R.~J. Stroeker and N.~Tzanakis.
\newblock Solving elliptic {D}iophantine equations by estimating linear forms
  in elliptic logarithms.
\newblock {\em Acta Arith.}, 67(2):177--196, 1994.

\bibitem{Tunnell}
J.~B. Tunnell.
\newblock A classical {D}iophantine problem and modular forms of weight
  {$3/2$}.
\newblock {\em Invent. Math.}, 72(2):323--334, 1983.

\bibitem{Tzanakis}
N.~Tzanakis.
\newblock Effective solution of two simultaneous {P}ell equations by the
  elliptic logarithm method.
\newblock {\em Acta Arith.}, 103(2):119--135, 2002.

\bibitem{Tz}
Nikos Tzanakis.
\newblock {\em Elliptic {D}iophantine equations}, volume~2 of {\em De Gruyter
  Series in Discrete Mathematics and Applications}.
\newblock Walter de Gruyter GmbH \& Co. KG, Berlin, 2013.
\newblock A concrete approach via the elliptic logarithm.

\end{thebibliography}

\end{document}